\documentclass[12pt]{amsart}
\usepackage{amssymb}
\usepackage{amsmath}
\usepackage{amscd}
\theoremstyle{plain}
\newtheorem{theorem}{Theorem}[section]
\newtheorem{corollary}[theorem]{Corollary}
\newtheorem{proposition}[theorem]{Proposition}
\newtheorem{lemma}[theorem]{Lemma}
{\theoremstyle{remark}

\newtheorem{remark}[theorem]{Remark}}
{\theoremstyle{definition}
\newtheorem{definition}[theorem]{Definition}
\newtheorem{example}[theorem]{Example}
}

\newcommand{\N}{\mathbb{N}}
\newcommand{\Z}{\mathbb{Z}}
\newcommand{\C}{\mathbb{C}}
\newcommand{\e}{\varepsilon}
\newcommand{\cK}{{\mathcal K}}

\newcommand{\M}{\mathbb{M}}
\newcommand{\F}{{\mathcal F}}
\newcommand{\Ca}{$C^*$-al\-ge\-bra }
\newcommand{\CA}{$C^*$-al\-ge\-bra}
\newcommand{\Csa}{$C^*$-sub\-al\-ge\-bra }
\newcommand{\CsA}{$C^*$-sub\-al\-ge\-bra}
\newcommand{\shom}{$*$-ho\-mo\-mor\-phism }

\begin{document}
\title{Non-separable AF-algebras}
\author{Takeshi Katsura}
\address{Department of Mathematics, 
Hokkaido University, Kita 10, Nishi 8, 
Kita-Ku, Sapporo, 060-0810, JAPAN}
\email{katsura@math.sci.hokudai.ac.jp}



\maketitle

\begin{abstract}
We give two pathological phenomena for non-separable AF-algebras 
which do not occur for separable AF-algebras. 
One is that non-separable AF-algebras are not determined 
by their Bratteli diagrams, 
and the other is that there exists a non-separable AF-algebra 
which is prime but not primitive. 
\end{abstract}

\section{Introduction}

In this paper, 
an AF-algebra means a \Ca 
which is an inductive limit of finite dimensional \CA s 
on {\em any} directed set. 
Equivalently, 

\begin{definition}
A \Ca $A$ is called an {\em AF-algebra} 
if it has a directed family of finite dimensional \CsA s 
whose union is dense in $A$. 
\end{definition}

When an AF-algebra $A$ is separable, 
we can find an increasing {\em sequence} of finite dimensional \CsA s 
whose union is dense in $A$. 
Thus for separable \CA s, 
the above definition coincides with the one in many literatures 
(for example, \cite{E}). 
For separable \CA s, 
there exists one more equivalent definition 
of AF-algebras: 

\begin{proposition}[{Theorem 2.2 of \cite{B}}]
A separable \Ca $A$ is an AF-algebra 
if and only if it is a locally finite dimensional \CA , 
which means that for any finite subset $\F$ of $A$ and any $\e>0$, 
we can find a finite dimensional \Csa $B$ of $A$ 
such that ${\rm dist}(x,B)<\e$ for all $x\in\F$. 
\end{proposition}

To the best of the author's knowledge, 
it is still open that the above lemma is valid in general. 

For each positive integer $n\in\Z_+$, 
$\M_n$ denotes the \Ca of all $n\times n$ matrices. 
Any finite dimensional \Ca $A$ 
is isomorphic to $\bigoplus_{i=1}^{k}\M_{n_i}$ 
for some $k\in\Z_+$ and ${}^{t}(n_1,\ldots,n_k)\in \Z_+^k$. 
Let $B\cong \bigoplus_{j=1}^{k'}\M_{n'_j}$ 
be another finite dimensional \CA . 
A \shom $\varphi\colon A\to B$ is determined up to unitary equivalence 
by the $k'\times k$ matrix $N$ 
whose $(j,i)$-entry is the multiplicity of 
the composition of the restriction of $\varphi$ 
to $\M_{n_i}\subset A_\lambda$ 
and the natural surjection from $B$ 
to $\M_{n'_j}$.

\begin{definition}
Let $\Lambda$ be a directed set with an order $\preceq$. 
An {\em inductive system of finite dimensional \CA s} 
$(A_\lambda,\varphi_{\mu,\lambda})$ over $\Lambda$ 
consists of a finite dimensional \Ca $A_\lambda$ 
for each $\lambda\in\Lambda$, 
and a \shom $\varphi_{\mu,\lambda}\colon A_\lambda\to A_{\mu}$ 
for each $\lambda,\mu\in\Lambda$ with $\lambda\prec\mu$ 
such that $\varphi_{\nu,\mu}\circ \varphi_{\mu,\lambda}
=\varphi_{\nu,\lambda}$ 
for $\lambda\prec\mu\prec\nu$. 

A {\em Bratteli diagram} of $(A_\lambda,\varphi_{\mu,\lambda})$ 
is the system $(n_\lambda,N_{\mu,\lambda})$ 
where 
$$
n_\lambda={}^{t}((n_\lambda)_1,\ldots, (n_\lambda)_{k_\lambda})
\in \Z_+^{k_\lambda}
$$ 
satisfies $A_\lambda\cong \bigoplus_{i=1}^{k_\lambda}\M_{(n_\lambda)_i}$ 
and $N_{\mu,\lambda}$ is 
$k_{\mu}\times k_\lambda$ matrix 
which indicates 
the multiplicities of the restrictions of $\varphi_{\mu,\lambda}$ 
as above. 
\end{definition}

A Bratteli diagram $(n_\lambda,N_{\mu,\lambda})$ satisfies 
$N_{\mu,\lambda}n_\lambda\leq n_\mu$ for $\lambda\prec\mu$, 
and $N_{\nu,\mu}N_{\mu,\lambda}=N_{\nu,\lambda}$ 
for $\lambda\prec\mu\prec\nu$. 
It is not difficult to see that 
when the directed set $\Lambda$ is $\Z_+$, 
any system $(n_\lambda,N_{\mu,\lambda})$ 
satisfying these two conditions can be 
realized as a Bratteli diagram of some 
inductive system of finite dimensional \CA s 
(see 1.8 of \cite{B}). 
This does not hold for general directed set: 

\begin{example}
Let $\Lambda=\{a,b,c,d,e\}$ 
with an order $a\succ b,c\succ d,e$. 
Let us define 
\[
n_a=(24), \quad 
n_b=\left(
\begin{matrix}
4\\ 4
\end{matrix}
\right),\quad 
n_c=\left(
\begin{matrix}
6\\ 6
\end{matrix}
\right),\quad 
n_d=\left(
\begin{matrix}
1\\ 3
\end{matrix}
\right),\quad 
n_e=\left(
\begin{matrix}
2\\ 2
\end{matrix}
\right),
\]
and 
\begin{align*}
N_{a,b}&=\left(
\begin{matrix}
3&3
\end{matrix}
\right),\qquad 
N_{b,d}=\left(
\begin{matrix}
1&1\\1&1
\end{matrix}
\right),\qquad 
N_{b,e}=\left(
\begin{matrix}
1&1\\1&1
\end{matrix}
\right),\qquad 
N_{a,d}=\left(
\begin{matrix}
6&6
\end{matrix}
\right),\\
N_{a,c}&=\left(
\begin{matrix}
2&2
\end{matrix}
\right),\qquad 
N_{c,d}=\left(
\begin{matrix}
3&1\\0&2
\end{matrix}
\right),\qquad 
N_{c,e}=\left(
\begin{matrix}
1&2\\2&1
\end{matrix}
\right),\qquad 
N_{a,e}=\left(
\begin{matrix}
6&6
\end{matrix}
\right). 
\end{align*}
These matrices satisfy $N_{\mu,\lambda}n_\lambda=n_\mu$ 
for $\lambda,\mu\in\Lambda$ with $\mu\succ \lambda$, 
and 
\[
N_{a,b} N_{b,d}=N_{a,c} N_{c,d}=N_{a,d},\qquad 
N_{a,b} N_{b,e}=N_{a,c} N_{c,e}=N_{a,e}.
\]
Thus the system 
$(n_\lambda, N_{\mu,\lambda})$ 
satisfies the two conditions above. 
However, one can see that 
this diagram never be a Bratteli diagram 
of inductive systems of finite dimensional \CA s. 
\end{example}

In 1.8 of \cite{B}, 
O. Bratteli showed that 
when the directed set $\Lambda$ is $\Z_+$, 
a Bratteli diagram of 
an inductive system of finite dimensional \CA s 
determines the inductive limit up to isomorphism. 
This is no longer true 
for general directed set $\Lambda$ 
as the following easy example shows. 

\begin{example}\label{example}
Let $X$ be an infinite set, 
and $\Lambda$ be the directed set consisting of 
all finite subsets of $X$ with inclusion as an order. 
We consider the following two 
inductive systems of finite dimensional \CA s. 

For each $\lambda\in \Lambda$, 
we define a \Ca 
$A_\lambda=\cK(\ell^2(\lambda))\cong \M_{|\lambda|}$ 
whose matrix unit is given by $\{e_{x,y}\}_{x,y\in \lambda}$. 
For $\lambda,\mu \in \Lambda$ with $\lambda\subset\mu$, 
we define a \shom $\varphi_{\mu,\lambda}\colon A_\lambda\to A_\mu$ 
by $\varphi_{\mu,\lambda}(e_{x,y})=e_{x,y}$. 
It is clear to see that this defines 
an inductive system of finite dimensional \CA s, 
and the inductive limit is $\cK(\ell^2(X))$. 

For each $\lambda\in \Lambda$ with $n=|\lambda|$, 
we set $A'_\lambda=\M_{n}$ 
whose matrix unit is given by $\{e_{k,l}\}_{1\leq k,l\leq n}$. 
For $\lambda,\mu \in \Lambda$ with $\lambda\subset\mu$, 
we define a \shom $\varphi'_{\mu,\lambda}\colon A'_\lambda\to A'_\mu$ 
by $\varphi'_{\mu,\lambda}(e_{k,l})=e_{k,l}$. 
It is clear to see that this defines 
an inductive system of finite dimensional \CA s, 
and the inductive limit is $\cK(\ell^2(\Z_+))$. 

The above two inductive systems give isomorphic 
Bratteli diagrams, 
but the AF-algebras $\cK(\ell^2(X))$ 
and $\cK(\ell^2(\Z_+))$ 
determined by the two inductive systems 
are isomorphic only when $X$ is countable. 

In a similar way, 
we can find two 
inductive systems of finite dimensional \CA s 
whose Bratteli diagrams are isomorphic, 
but the inductive limits are $\bigotimes_{x\in X}\M_2$ 
and $\bigotimes_{k=1}^{\infty}\M_2$ 
which are not isomorphic when $X$ is uncountable. 
\end{example}

By Example \ref{example}, 
we can see that 
G. A. Elliott's celebrated theorem 
of classifying (separable) AF-algebras using $K_0$-groups 
(Theorem 6.4 of \cite{E}) does not follow 
for non-separable AF-algebras, 
because $K_0$-groups are determined by Bratteli diagrams. 
Example \ref{example} is not so interesting 
because the inductive system $(A'_\lambda,\varphi'_{\mu,\lambda})$ 
has many redundancies and does not come from 
directed families of finite dimensional \CsA s. 
More interestingly, 
we can get the following 
whose proof can be found in the next section: 

\begin{theorem}\label{main1}
There exist two non-isomorphic AF-algebras $A$ and $B$ 
such that 
they have directed families of finite dimensional \CsA s 
which define isomorphic Bratteli diagrams. 
\end{theorem}

The author could not find such an example 
in which every finite dimensional \CsA s 
are isomorphic to full matrix algebras $\M_n$ 
(cf.\ Problem 8.1 of \cite{D}). 

As another pathological fact on 
non-separable AF-algebras, 
we prove the next theorem in Section \ref{Sec:PP}.

\begin{theorem}\label{main2}
There exists a non-separable AF-algebra 
which is prime but not primitive.
\end{theorem}

It had been a long standing problem 
whether there exists a \Ca which is prime but not primitive, 
until N. Weaver found such a \Ca in \cite{W}. 
Note that such a \Ca cannot be separable. 

\medskip

{\bf Acknowledgments.} 
The author is grateful 
to the organizers of the Abel Symposium 2004 
for giving him opportunities to talk in the conference 
and to contribute in this volume. 
He is also grateful to George A. Elliott and Akitaka Kishimoto 
for useful comments. 
This work was partially supported by Research Fellowship 
for Young Scientists of the Japan Society for the Promotion of Science.

\section{Proof of Theorem \ref{main1}}

In this section, 
we will prove Theorem \ref{main1}. 
Let $X$ be an infinite set, 
and $Z$ be the set of all subsets $z$ of $X$ 
with $|z|=2$. 

For each $z\in Z$, 
we define a \Ca $M_z$ by $M_z=\M_2$. 
Elements of the direct product $\prod_{z\in Z}M_z$ 
will be considered as norm bounded functions $f$ on $Z$ 
such that $f(z)\in M_z$ for $z\in Z$. 
For each $z\in Z$, 
we consider $M_z\subset \prod_{z\in Z}M_z$ 
as a direct summand. 
We denote by $\bigoplus_{z\in Z}M_z$ 
the direct sum of $M_z$'s 
which is an ideal of $\prod_{z\in Z}M_z$. 

\begin{definition}
For each $z\in Z$, 
we fix a matrix unit $\{e_{i,j}^{z}\}_{i,j=1}^2$ of $M_z=\M_2$. 
For each $x\in X$, 
we define a projection $p_x\in \prod_{z\in Z}M_z$ 
by 
$$p_x(z)=\begin{cases}
e_{1,1}^{z} & \text{if $x\in z$,}\\
0 & \text{if $x\notin z$.}
\end{cases}$$
We denote by $A$ the \Csa of $\prod_{z\in Z}M_z$ 
generated by $\bigoplus_{z\in Z}M_z$ and $\{p_x\}_{x\in X}$. 
\end{definition}

\begin{definition}
For each $z=\{x_1,x_2\}\in Z$, 
we fix a matrix unit $\{e_{x_i,x_j}^{z}\}_{i,j=1}^2$ 
of $M_z=\M_2$. 
For each $x\in X$, 
we define a projection $q_x\in \prod_{z\in Z}M_z$ 
by 
$$q_x(z)=\begin{cases}
e_{x,x}^{z} & \text{if $x\in z$,}\\
0 & \text{if $x\notin z$.}
\end{cases}$$
We denote by $B$ the \Csa of $\prod_{z\in Z}M_z$ 
generated by $\bigoplus_{z\in Z}M_z$ and $\{q_x\}_{x\in X}$. 
\end{definition}

The following easy lemma illustrates an difference of $A$ and $B$. 

\begin{lemma}
For $x,y\in X$ with $x\neq y$, 
we have $p_xp_y=e_{1,1}^{\{x,y\}}\neq 0$, 
and $q_xq_y=0$. 
\end{lemma}

\begin{proof}
Straightforward. 
\end{proof}

\begin{definition}
Let $\lambda$ be a finite subset of $X$. 
We denote by $A_{\lambda}$ the \Csa of $A$ spanned by 
$\bigoplus_{z\subset\lambda}M_z$ and $\{p_x\}_{x\in\lambda}$, 
and by $B_{\lambda}$ the \Csa of $B$ spanned by 
$\bigoplus_{z\subset\lambda}M_z$ and $\{q_x\}_{x\in\lambda}$, 
\end{definition}

\begin{lemma}\label{Lem:finite}
There exist isomorphisms 
\[
A_{\lambda}\cong B_{\lambda}\cong 
\bigoplus_{z\subset\lambda}\M_2\oplus \bigoplus_{x\in\lambda}\C
\] 
for each finite set $\lambda\subset X$ 
such that two inclusions $A_{\lambda}\subset A_{\mu}$ 
and $B_{\lambda}\subset B_{\mu}$ 
have the same multiplicity. 
\end{lemma}

\begin{proof}
For $x\in\lambda$, 
let us denote $p'_x\in A_{\lambda}$ 
by $p'_x=p_x-\sum_{y\in \lambda\setminus\{x\}}e_{1,1}^{\{x,y\}}$. 
Then we have an orthogonal decomposition 
$$A_{\lambda}=\sum_{z\subset\lambda}M_z + \sum_{x\in\lambda}\C p'_x.$$
This proves $A_{\lambda}\cong
\bigoplus_{z\subset\lambda}\M_2\oplus \bigoplus_{x\in\lambda}\C$. 
Similarly we have $B_{\lambda}\cong
\bigoplus_{z\subset\lambda}\M_2\oplus \bigoplus_{x\in\lambda}\C$. 
Now it is routine to check the last statement. 
\end{proof}

\begin{proposition}\label{Prop:sameBD}
Two \CA s $A$ and $B$ are AF-algebras, 
and the directed families $\{A_\lambda\}$ 
and $\{B_\lambda\}$ of finite dimensional \CsA s 
give isomorphic Bratteli diagrams. 
\end{proposition}

\begin{proof}
Follows from the facts 
$$\overline{\bigcup_{\lambda\subset X}A_\lambda}=A,\quad 
\overline{\bigcup_{\lambda\subset X}B_\lambda}=B$$
and Lemma \ref{Lem:finite}.
\end{proof}

\begin{remark}
From Proposition \ref{Prop:sameBD}, 
we can show that $K_0(A)$ and $K_0(B)$ 
are isomorphic as scaled ordered groups.
In fact, they are isomorphic to 
the subgroup $G$ of $\prod_{z\in Z}\Z$ 
generated by $\bigoplus_{z\in Z}\Z$ and $\{g_x\}_{x\in X}$, 
where $g_x\in \prod_{z\in Z}\Z$ is defined by 
\[
g_x(z)=\begin{cases}
1 & \text{if $x\in z$,}\\
0 & \text{if $x\notin z$.}
\end{cases}
\]
The order of $G$ is the natural one, and its scale is 
\[
\{g\in G\mid \text{$0\leq g(z)\leq 2$ 
for all $z\in Z$}\}.
\]
From this fact and Elliott's theorem (Theorem 6.4 of \cite{E}), 
we can show the next lemma, 
although we give a direct proof here. 
\end{remark}

\begin{proposition}\label{Prop:count}
When $X$ is countable, 
$A$ and $B$ are isomorphic. 
\end{proposition}

\begin{proof}
Let us list $X=\{x_1,x_2,\ldots\}$. 
We define a \shom $\varphi\colon A\to B$ as follows. 
For $z=\{x_k,x_l\}$, 
we define $\varphi(e_{i,j}^{z})=e_{x_{n_i},x_{n_j}}^{z}$ 
where $n_1=k,n_2=l$ when $k<l$ and $n_1=l,n_2=k$ when $k>l$. 
For $x_k\in X$, we set 
\[
\varphi(p_{x_k})=q_{x_k}+\sum_{i=1}^{k-1}
\big(e_{x_{i},x_{i}}^{\{x_i,x_k\}}-e_{x_{k},x_{k}}^{\{x_i,x_k\}}\big).
\]
Now it is routine to check that $\varphi$ is an isomorphism 
from $A$ to $B$. 
\end{proof}

Proposition \ref{Prop:count} is no longer true for uncountable $X$. 
To see this, we need the following lemma. 

\begin{lemma}\label{lem:exact}
There exists a surjection 
$\pi_A\colon A\to \bigoplus_{x\in X}\C$ 
defined by $\pi_A(M_z)=0$ for $z\in Z$ 
and $\pi_A(p_x)=\delta_x$ for $x\in X$. 
Its kernal is $\bigoplus_{z\in Z}M_z$ 
which coincides with 
the ideal generated 
by the all commutators $xy-yx$ of $A$.  
The same is true for $B$. 
\end{lemma}

\begin{proof}
Let $\pi_A$ be the quotient map from $A$ to 
$A/\bigoplus_{z\in Z}M_z$. 
Then $A/\bigoplus_{z\in Z}M_z$ 
is generated by $\{\pi_A(p_x)\}_{x\in X}$ 
which is an orthogonal family of non-zero projections. 
This proves the first statement. 
Since $\bigoplus_{x\in X}\C$ is commutative, 
the ideal $\bigoplus_{z\in Z}M_z$ 
contains all commutators. 
Conversely, the ideal generated by the commutators of $A$ 
contains $\bigoplus_{z\in Z}M_z$ 
because $\M_2$ is generated by its commutators. 
This shows that $\bigoplus_{z\in Z}M_z$ 
is the ideal generated by the all commutators of $A$. 
The proof goes similarly for $B$. 
\end{proof}

\begin{proposition}\label{Prop:non-isom}
When $X$ is uncountable, 
$A$ and $B$ are not isomorphic. 
\end{proposition}

\begin{proof}
To the contrary, 
suppose that there exists an isomorphism $\varphi\colon A\to B$. 
By Lemma \ref{lem:exact}, 
$\bigoplus_{z\in Z}M_z$ 
is the ideal generated by the all commutators in both $A$ and $B$. 
Hence $\varphi$ preserves this ideal 
$\bigoplus_{z\in Z}M_z$. 
Thus we get the following commutative diagram with exact rows; 
\[
\begin{CD}
0 @>>> \bigoplus_{z\in Z}M_z 
@>>> A @>\pi_A>> \bigoplus_{x\in X}\C @>>> 0\phantom{.} \\
@. @VV\varphi V @VV\varphi V  @VVV \\
0 @>>> \bigoplus_{z\in Z}M_z 
@>>> B @>\pi_B>> \bigoplus_{x\in X}\C @>>> 0. \\
\end{CD}
\]
Since the family of projections $\{q_x\}_{x\in X}$ in $B$ 
is mutually orthogonal, 
the surjection $\pi_B\colon B\to \bigoplus_{x\in X}\C$ 
has a splitting map $\sigma_B\colon \bigoplus_{x\in X}\C\to B$ 
defined by $\sigma_B(\delta_x)=q_x$. 
Hence by the diagram above, 
the surjection $\pi_A\colon A\to \bigoplus_{x\in X}\C$ 
also has a splitting map $\sigma_A\colon \bigoplus_{x\in X}\C\to A$. 
Let us set $p_x'=\sigma_A(\delta_x)$ for $x\in X$. 
Choose a countable infinite subset $Y$ of $X$. 
For each $y\in Y$, 
the set 
\[
\F_y=\big\{x\in X\ \big|\ x\neq y, \|(p_y-p_y')(\{x,y\})\|\geq 1/2\big\}
\]
is finite, 
because $p_y-p_y'\in \ker\pi_A=\bigoplus_{z\in Z}M_z$. 
Since $X$ is uncountable, 
we can find $x_0\in X$ with $x_0\notin Y\cup\bigcup_{y\in Y}\F_y$. 
Since 
\[
\F_{x_0}=
\big\{x\in X\ \big|\ x\neq x_0, \|(p_{x_0}-p_{x_0}')(\{x,x_0\})\|\geq 1/2\big\}
\]
is finite, 
we can find $y_0\in Y\setminus \F_{x_0}$. 
We set $z=\{x_0,y_0\}$. 
From $y_0\notin\F_{x_0}$, 
we have $\|(p_{x_0}-p_{x_0}')(z)\|< 1/2$, 
and from $x_0\notin\F_{y_0}$, 
we have $\|(p_{y_0}-p_{y_0}')(z)\|< 1/2$. 
However, $p_{x_0}(z)=p_{y_0}(z)=e_{1,1}^{z}$ 
and $p_{x_0}'(z)$ is orthogonal to $p_{y_0}'(z)$. 
This is a contradiction. 
Thus $A$ and $B$ are not isomorphic. 
\end{proof}

Combining Proposition \ref{Prop:sameBD} and 
Proposition \ref{Prop:non-isom}, 
we get Theorem \ref{main1}.

\section{A prime AF-algebra which is not primitive}\label{Sec:PP}

In this section, 
we construct an AF-algebra which is prime 
but not primitive. 
Although we follow the idea of Weaver in \cite{W}, 
our construction of the \Ca 
and proof of the main theorem is much easier 
than the ones there. 
A similar construction can be found in \cite{K}, 
but the proof there uses general facts 
of topological graph algebras.

Let $X$ be an uncountable set, 
and $\Lambda$ be the directed set of 
all finite subsets of $X$. 
For $n\in\N$, we set 
$\Lambda_n=\big\{\lambda\subset X\ \big|\  |\lambda|=n\big\}$.
We get $\Lambda=\coprod_{n=0}^\infty \Lambda_n$.

\begin{definition}
For $n\in\Z_+$ and $\lambda\in \Lambda_n$, 
we define 
\[
l(\lambda)=\{t\colon\{1,\ldots,n\}\to\lambda\mid 
\text{$t$ is a bijection}\}.
\]
For $\emptyset\in\Lambda$, 
we define $l(\emptyset)=\{\emptyset\}$. 
\end{definition}

Note that $|l(\lambda)|=n!$ for $\lambda\in \Lambda_n$ and $n\in\N$. 

\begin{definition}
For $n\in\N$ and $\lambda\in \Lambda_n$, 
we define $M_\lambda\cong \M_{n!}$ 
whose matrix unit is given by 
$\{e_{s,t}^{(\lambda)}\}_{s,t\in l(\lambda)}$. 
\end{definition}

\begin{definition}
Take $\lambda\in \Lambda_n$ and $\mu\in \Lambda_m$ with $\lambda\cap\mu=\emptyset$. 
For $t\in l(\lambda)$ and $s\in l(\mu)$, 
we define $ts\in l(\lambda\cup\mu)$ by 
\[
(ts)(i)=\begin{cases}
t(i)&\text{for $i=1,\ldots,n$}\\
s(i-n)&\text{for $i=n+1,\ldots,n+m$}.
\end{cases}
\]
\end{definition}

Note that when $\mu=\emptyset$, we have $t\emptyset=t$. 

\begin{definition}
For $\lambda,\mu\in\Lambda$ with $\lambda\subset\mu$, 
we define a \shom $\iota_{\mu,\lambda}\colon M_\lambda\to M_\mu$ 
by
\[
\iota_{\mu,\lambda}\big(e_{s,t}^{(\lambda)}\big)
=\sum_{u\in l(\mu\setminus\lambda)}e_{su,tu}^{(\mu)}
\quad \text{for $s,t\in l(\lambda)$.}
\]
\end{definition}

Note that $\iota_{\lambda,\lambda}$ is 
the identity map of $M_\lambda$, 
and that 
$\iota_{\lambda_3,\lambda_2}\circ\iota_{\lambda_2,\lambda_1}\neq 
\iota_{\lambda_3,\lambda_1}$ 
for $\lambda_1\subsetneq\lambda_2\subsetneq\lambda_3$. 
For $\lambda_1,\lambda_2\in \Lambda_n$ and $\mu\in \Lambda_m$ 
with $\lambda_1\neq \lambda_2$ and $\lambda_1\cup\lambda_2\subset\mu$, 
the images $\iota_{\mu,\lambda_1}(M_{\lambda_1})$ 
and $\iota_{\mu,\lambda_2}(M_{\lambda_2})$ 
are mutually orthogonal. 

\begin{definition}
For $\lambda\in\Lambda$, 
we define a \shom 
$\iota_\lambda\colon M_\lambda\to \prod_{\mu\in\Lambda}M_\mu$ 
by 
\[
\iota_\lambda(x)(\mu)=\begin{cases}
\iota_{\mu,\lambda}(x) & \text{if $\lambda\subset \mu$,}\\
0 & \text{otherwise,}
\end{cases}
\]
for $x\in M_\lambda$. 
We set $N_\lambda=\iota_\lambda(M_\lambda)\subset \prod_{\mu\in\Lambda}M_\mu$ 
and $f_{s,t}^{(\lambda)}=\iota_\lambda(e_{s,t}^{(\lambda)})\in N_\lambda$ 
for $s,t\in l(\lambda)$. 
\end{definition}

For $\lambda\in\Lambda_n$, 
We have $N_\lambda\cong \M_{n!}$ 
and $\{f_{s,t}^{(\lambda)}\}_{s,t\in l(\lambda)}$ 
is a matrix unit of $N_\lambda$. 

\begin{lemma}\label{easy1}
For $\lambda,\mu\in\Lambda$ with $\lambda\subset\mu$, 
$s,t\in l(\lambda)$ and $s',t'\in l(\mu)$, 
we have 
$f_{s,t}^{(\lambda)}f_{s',t'}^{(\mu)}=f_{su,t'}^{(\mu)}$ 
when $s'=tu$ with some $u\in l(\mu\setminus \lambda)$, 
and $f_{s,t}^{(\lambda)}f_{s',t'}^{(\mu)}=0$ otherwise. 
\end{lemma}

\begin{proof}
Straightforward. 
\end{proof}

\begin{lemma}\label{easy2}
For $\lambda,\mu\in\Lambda$, 
we have $0\neq N_\lambda N_{\mu}\subset N_{\mu}$ 
if $\lambda\subset\mu$, 
$0\neq N_\lambda N_{\mu}\subset N_{\lambda}$ 
if $\lambda\supset\mu$, 
and $N_\lambda N_{\mu}=0$ otherwise. 
\end{lemma}

\begin{proof}
If $\lambda\subset\mu$, 
we have $0\neq N_\lambda N_{\mu}\subset N_{\mu}$ 
by Lemma \ref{easy1}. 
Similarly we have 
$0\neq N_\lambda N_{\mu}\subset N_{\lambda}$ 
if $\lambda\supset\mu$. 
Otherwise, 
we can easily see $N_\lambda N_{\mu}=0$ 
from the definition. 
\end{proof}

\begin{corollary}\label{easy3}
For each $n$, 
the family $\{N_\lambda\}_{\lambda\in \Lambda_n}$ of \CA s
is mutually orthogonal. 
\end{corollary}

\begin{corollary}\label{easy4}
Take $\lambda,\lambda'\in\Lambda$ with $\lambda\subset\lambda'$. 
Let $p_{\lambda'}$ be the unit of $N_{\lambda'}$. 
Then $N_{\lambda}\ni a\mapsto ap_{\lambda'}\in N_{\lambda'}$ 
is an injective $*$-ho\-mo\-mor\-phism. 
\end{corollary}

\begin{definition}
We define 
$A=\overline{\sum_{\lambda\in\Lambda}N_\lambda}
\subset \prod_{\mu\in\Lambda}M_\mu$. 
\end{definition}

\begin{proposition}\label{AF}
The set $A$ is an AF-algebra. 
\end{proposition}

\begin{proof}
For each $\mu\in\Lambda$, 
$A_\mu=\sum_{\lambda\subset\mu}N_\lambda$ 
is a finite dimensional \Ca by Lemma \ref{easy2}. 
For $\lambda,\mu\in\Lambda$ with $\lambda\subset \mu$, 
we have $A_\lambda\subset A_\mu$. 
Hence $A=\overline{\bigcup_{\mu\in\Lambda}A_\mu}$ 
is an AF-algebra. 
\end{proof}

\begin{lemma}\label{non-zero}
Every non-zero ideal $I$ of $A$ 
contains $N_\lambda$ for some $\lambda\in\Lambda$. 
\end{lemma}

\begin{proof}
As in the proof of Proposition \ref{AF}, 
we set $A_\mu=\sum_{\lambda\subset\mu}N_\lambda$ for $\mu\in\Lambda$. 
Since $A=\overline{\bigcup_{\mu\in\Lambda}A_\mu}$, 
we have $I=\overline{\bigcup_{\mu\in\Lambda}(I\cap A_\mu)}$ 
for an ideal $I$ of $A$. 
Hence if $I$ is nonzero, 
we have $I\cap A_{\mu_0}\neq 0$ for some $\mu_0\in\Lambda$. 
Thus we can find a non-zero element $a\in I$ 
in the form $a=\sum_{\lambda\subset \mu_0}a_\lambda$ 
for $a_\lambda\in N_\lambda$. 
Since $a\neq 0$, we can find $\lambda_0\in \Lambda$ 
with $\lambda_0\subset \mu_0$ 
such that $a_{\lambda_0}\neq 0$ 
and $a_\lambda=0$ for all $\lambda\subsetneq\lambda_0$. 
Take $x_0\in X$ with $x_0\notin\mu_0$. 
Set $\lambda_0'=\lambda_0\cup\{x_0\}$. 
Let $p_{\lambda_0'}$ be the unit of $N_{\lambda_0'}$. 
For $\lambda\subset \mu_0$, 
$a_\lambda p_{\lambda_0'}\neq 0$ 
only when $\lambda\subset \lambda_0$. 
Hence we have $ap_{\lambda_0'}=a_{\lambda_0}p_{\lambda_0'}$. 
By Corollary \ref{easy4}, 
$a_{\lambda_0}p_{\lambda_0'}$ 
is a non-zero element of $N_{\lambda_0'}$. 
Hence we can find a non-zero element in $I\cap N_{\lambda_0'}$. 
Since $N_{\lambda_0'}$ is simple, 
we have $N_{\lambda_0'}\subset I$. 
We are done. 
\end{proof}

\begin{lemma}\label{hereditary}
If an ideal $I$ of $A$ satisfies $N_{\lambda_0}\subset I$ 
for some $\lambda_0\in\Lambda$, 
then $N_{\lambda}\subset I$ for all $\lambda\supset \lambda_0$. 
\end{lemma}

\begin{proof}
Clear from Lemma \ref{easy2} and the simplicity of $N_{\lambda}$. 
\end{proof}

\begin{proposition}
The \Ca is prime but not primitive. 
\end{proposition}

\begin{proof}
Take two non-zero ideals $I_1,I_2$ of $A$. 
By Lemma \ref{non-zero}, 
we can find $\lambda_1,\lambda_2\in\Lambda$ 
such that $N_{\lambda_1}\subset I_1$ 
and $N_{\lambda_2}\subset I_2$. 
Set $\lambda=\lambda_1\cup\lambda_2\in \Lambda$. 
By Lemma \ref{hereditary}, 
we have $N_{\lambda}\subset I_1\cap I_2$. 
Thus $I_1\cap I_2\neq 0$. 
This shows that $A$ is prime. 

To prove that $A$ is not primitive, 
it suffices to see that for any state $\varphi$ of $A$ 
we can find a non-zero ideal $I$ such that $\varphi(I)=0$ 
(see \cite{W}). 
Take a state $\varphi$ of $A$. 
By Corollary \ref{easy3}, 
the family $\{N_\lambda\}_{\lambda\in \Lambda_n}$ of \CA s
is mutually orthogonal for each $n\in\N$. 
Hence the set 
\[
\Omega_n=\{\lambda\in \Lambda_n\mid 
\text{the restriction of $\varphi$ to $N_{\lambda}$ is non-zero}\}
\]
is countable for each $n\in\N$. 
Since $X$ is uncountable, 
we can find $x_0\in X$ 
such that $x_0\notin\lambda$ 
for all $\lambda\in\bigcup_{n\in\N}\Omega_n$. 
Let $I=\overline{\sum_{\lambda\ni x_0}N_\lambda}$. 
Then $I$ is an ideal of $A$ by Lemma \ref{easy2}. 
Since $\lambda\ni x_0$ implies $\varphi(N_\lambda)=0$, 
we have $\varphi(I)=0$. 
Therefore $A$ is not primitive. 
\end{proof}

This finishes the proof of Theorem \ref{main2}. 

\begin{remark}
Let $(A_\lambda,\varphi_{\mu,\lambda})$ 
be an inductive system of finite dimensional \CA s 
over a directed set $\Lambda$, 
and $A$ be its inductive limit. 
It is not hard to see that 
the AF-algebra $A$ is prime if and only if 
the Bratteli diagram of the inductive system 
satisfies the analogous condition of (iii) 
in Corollary 3.9 of \cite{B}. 
Hence, 
the Bratteli diagram of 
an inductive system of finite dimensional \CA s 
determines the primeness of the inductive limit, 
although it does not determine the inductive limit itself. 
However the primitivity of the inductive limit 
is not determined by the Bratteli diagram. 
In fact, 
in a similar way to the construction of Example 
\ref{example}, 
we can find an inductive system of finite dimensional \CA s 
whose Bratteli diagram is isomorphic to the one 
coming from the directed family $\{A_\lambda\}$ 
constructed in the proof of Proposition \ref{AF}, 
but the inductive limit is separable. 
This AF-algebra is primitive 
because it is separable and prime 
(see, for example, Proposition 4.3.6 of \cite{P}). 
\end{remark}


\begin{thebibliography}{[W03]}

\bibitem[B72]{B}
Bratteli, O. {\it Inductive limits of finite dimensional $C\sp{*} $-algebras.} Trans. Amer. Math. Soc. 171 (1972), 195--234. 

\bibitem[D67]{D}
Dixmier, J. {\it On some $C\sp{*} $-algebras considered by Glimm.} J. Funct. Anal. {\bf 1} (1967) 182--203.

\bibitem[E76]{E}
Elliott, G. A. {\it On the classification of inductive limits of sequences of semisimple finite-dimensional algebras.} J. Algebra {\bf 38} (1976), no. 1, 29--44.

\bibitem[K04]{K}
Katsura, T. {\it A class of $C^*$-algebras generalizing both graph algebras and homeomorphism $C^*$-algebras III, ideal structures.} Preprint 2004, math.OA/0408190. 

\bibitem[P79]{P}
Pedersen, G. K. {\it $C\sp{*} $-algebras and their automorphism groups.} London Mathematical Society Monographs, {\bf 14}. Academic Press, Inc., London-New York, 1979.

\bibitem[W03]{W}
Weaver, N. {\it A prime $C\sp *$-algebra that is not primitive.} J. Funct. Anal. {\bf 203} (2003), no. 2, 356--361.

\end{thebibliography}
\end{document}